\documentclass[11pt,a4paper]{article}
 \usepackage{euscript}

\usepackage{amsmath}
\usepackage{graphicx}
\usepackage{amsthm}
\usepackage{amssymb}
\usepackage{epsfig}
\usepackage[all]{xy}
\usepackage{url}

\theoremstyle{theorem}
\newtheorem{theorem}{Theorem}[section]
\newtheorem{corollary}[theorem]{Corollary}

\newtheorem{lemma}[theorem]{Lemma}
\newtheorem{proposition}[theorem]{Proposition}

\newtheorem{remark}{Remark}[section]
\numberwithin{equation}{section}

\title{The Penrose Transform in the Split Signature}
\date{\today}
\author{Masood Aryapoor\footnote{masood.aryapoor@yale.edu}\\ \\
Mathematics Department\\
Yale University\\
442 Dunham Lab\\
10 Hillhouse Avenue\\
New Haven, CT 06511 USA\\}

\begin{document}  
\maketitle
 \begin{abstract}
A version of the Penrose transform is introduced in the split signature. It relates the cohomological data on $\mathbb{CP}^3\setminus\mathbb{RP}^3$ and kernel of differential operators on $M$, the (real) Grassmannian of 2-planes in $\mathbb{R}^4$. As an example we derive the following cohomological interpretation of the so-called X-ray transform 
$$H^1_c(\mathbb{C}\mathbb{P}^3\setminus\mathbb{R}\mathbb{P}^3,\mathcal{O}(-2))\overset{\mathcal{\cong}}{\to}\ker\left (\square_{2,2}:\Gamma^{\omega}(
M,\widetilde{\varepsilon[-1]})\to
\Gamma^{\omega}(M,\widetilde{\varepsilon[-3]})\right )
$$
where 
$\Gamma^{\omega}(
M,\widetilde{\varepsilon[-1]})$ and $ \Gamma^{\omega}(M,\widetilde{\varepsilon[-3]})$ are real analytic sections of certain (homogeneous) line bundles on $M$, $c$ stands for cohomology with compact support and $\square_{2,2}$ is the ultrahyperbolic operator. Furthermore, this gives a cohomological realization of the so-called "minimal" representation of $SL(4,\mathbb{R})$. We also present the split Penrose transform in split instanton backgrounds.
\end{abstract}
\begin{section}{Introduction} 
The Penrose transform is a well-known transform  which relates   the cohomological data on $\mathbb{P}$, the complex 3-projective space,  and the spaces of solutions of certain differential operators on $\mathbb{M}$, the Grassmannian of complex 2-planes in $\mathbb{C}^4$, i.e. the complexification of the conformal compactification of the Minkowski space, see \cite{BE,EPW}.  On the other hand, there is a well-known transform in Real Integral Geometry which is called the X-ray (or Radon) transform: a smooth function $f$ (or section of an appropriate line bundle) on the totally real submanifold $P=\mathbb{RP}^3$ can be integrated along lines to yield a function $\phi$ on $M$, the Grassmannian of 2-planes in $\mathbb{R}^4$, see \cite{He,Sp}. It is a classical result that $\phi$ is a solution to the ultrahyperbolic wave equation and that all 
such solutions determine a unique $f$ on $\mathbb{RP}^3$, see \cite{Jo}. Following Atiyah, locally one can think of the function $f$ (when $f$ is real analytic) as a preferred Cech cocycle, see \cite{At}. However, globally this cohomological description  of the X-ray transform breaks down. Nevertheless, this suggests that there might be a relation between the X-ray transform and the Penrose transform. There has been  substantial research on finding the precise relationship between the Penrose transform and the X-ray transform, see \cite{Ea,GS,Wo,Sp,BEGM,BEGM1}.  As we will see in this paper, the cohomological interpretation of the X-ray transform via the Penrose transform is obtained  by working with  Cohomology Theory  with compact support rather than the usual Cohomology Theory, see transform \ref{CI}.\\
In \cite{Ea},  the authors have introduced a version of the Penrose transform which yields a family of transformations in Real Integral Geometry one of which is the X-ray transform. 
In this paper we introduce another version of the Penrose transform in the split signature which deals with cohomological data.
There are several advantages to this cohomological interpretation. First of all, the ingredients of this approach are almost as those in the Penrose transform and we do not need to deal with less known concepts such as "involutive" structures as in \cite{Ea}.  It is relatively simple because a major part of it can be done using the complex Penrose transform with no extra work. It makes it possible to define the X-ray transform locally as in the Penrose transform. Moreover, the data in the split Penrose transform is in parallel with the one in the usual Penrose transform.\\
Now we explain various real forms of the Penrose transform, see \cite{DM}. 
The group $SL(4,\mathbb{C})$ acts on the spaces involved in the Penrose transform and, furthermore, the Penrose transform is  $SL(4,\mathbb{C})$-equivariant.   It is interesting to look at the real forms of the Penrose transform. In the introduction, for simplicity, we only consider the Penrose transform for $\mathcal{O}(-2)$ which is a very important case.\\
\noindent There are three real forms:\\
 (1) \textbf{Euclidean}:  This corresponds to the real form $SO_0(5,1)$ of $SL(4,\mathbb{C})$ which gives the totally real sub-manifold $S^4$ of $\mathbb{M}$.  
In this case we have the fibration  $\pi:\mathbb{P}\to S^4$, see \cite{At}.  Using this fibration, we have the following Penrose transform in the Euclidean case
$$H^1(\pi^{-1}(U),\mathcal{O}(-2))\to \ker \square|_{U}$$
which is an isomorphism. Here $\ker \square|_{U}$ is just the space of harmonic functions (more precisely sections of a line bundle) on $U$. Note that, globally we have  the isomorphism
$$H^1(\mathbb{P},\mathcal{O}(-2))\to \ker \square_{S^4}$$
which is trivial because both spaces are zero. In contrast, as we will see, the Penrose transform produces very interesting global isomorphisms in other real forms. For a detailed discussion of the Penrose transform in the Euclidean picture see \cite{At}. \\
 \noindent (2) \textbf{Minkowski}: This corresponds to the real form $SU(2,2)$ of $SL(4,\mathbb{C})$ which gives the totally real sub-manifold $M_0 \cong S^1\times S^3$, the compactified Minkowski space, of $\mathbb{M}$. There is a Penrose transform in this case which is due to Wells, see \cite{We}. It gives  a bijective Penrose transform
 $$H^1(Q,\mathcal{O}(-2))\to\ker \square_{3,1}^{\omega}$$
 where $Q$ is a certain five dimensional real submanifold of $\mathbb{P}$, $\square_{3,1}$ is just the wave operator on $M_0$ and ${\omega}$ stands for  real analytic solutions. In order to construct this transform, one can consider small neighborhoods of $Q$ in $\mathbb{P}$, apply the usual Penrose transform and then restrict the data to $M_0$. It is easy to see that this transform is injective, but
 the hard part is to prove that this transform is surjective.\\
\noindent (3) \textbf{Split}:  This corresponds to the real form $SL(4,\mathbb{R})$ of $SL(4,\mathbb{C})$ which gives the totally real sub-manifold $M=Gr(2,\mathbb{R}^4)$, the Grassmannian of 2-planes in $\mathbb{R}^4$, of $\mathbb{M}$.  This case was known in a particular form before the appearance of the Penrose transform. More precisely, there is a well-known transform
$$\Gamma(\mathbb{R}\mathbb{P}^3,\varepsilon(-2))\overset{\mathcal{R}}{\to} \ker{\square_{2,2}}$$
where $\Gamma(\mathbb{R}\mathbb{P}^3,\varepsilon(-2))$ is the space of smooth homogeneous functions on $\mathbb{R}^4\setminus 0$ of homogeneity -2 and $\square_{2,2}$ is the ultrahyperbolic operator on $M$ acting on appropriate line bundles on $M$, see section 6.
 This transform is known as the X-ray transform or Radon transform, see \cite{He,Sp}.  It is well-known that this transform gives a bijection. 
It is possible to prove this result using a version of the Penrose transform adapted for the split case, see \cite{Ea}. \\ 
 
\noindent In all these versions of the Penrose transform (complex, Euclidean, Minkowski and split), the output of the transform is the space of solutions of the corresponding Laplacian.  Moreover the input of the transformation is, roughly speaking, the first cohomology group of $\mathcal{O}(-2)$  except in the split case. It is natural to ask if there is a version of the Penrose transform in the split signature relating the first cohomology group of $\mathcal{O}(-2)$ to the space of solutions of the ultrahyperbolic equation.
Here we propose a version of the split Penrose transform which deals with the first cohomology group with supports. As an example we show that there is an isomorphism
\begin{equation}\label{CI}
H^1_c(\mathbb{C}\mathbb{P}^3\setminus\mathbb{R}\mathbb{P}^3,\mathcal{O}(-2))\overset{\mathcal{P}}{\to}\ker\left (\square_{2,2}:\Gamma^{\omega}(
M,\widetilde{\varepsilon[-1]})\to
\Gamma^{\omega}(M,\widetilde{\varepsilon[-3]})\right )
\end{equation}
where $\Gamma^{\omega}(
M,\widetilde{\varepsilon[-1]})$ and $ \Gamma^{\omega}(M,\widetilde{\varepsilon[-3]})$ are real analytic sections of certain (homogeneous) line bundles on $M$ and $c$ stands for cohomology with compact support.
It is worth noting that $H^1(\mathbb{CP}^3\setminus\mathbb{R}\mathbb{P}^3,\mathcal{O}(-2))$ is zero, see \cite{Ha}. \\
It is possible to obtain a local version of the split Penrose transform as well. 
More precisely we will show that for  any open  subset $U$ of $M$ there is a map
$$H^1_\Phi(U'',\mathcal{O}(-2))\overset{\mathcal{P}}{\to}\ker\left (\square_{2,2}:\Gamma^{\omega}(
U,\widetilde{\varepsilon[-1]})\to
\Gamma^{\omega}(U,\widetilde{\varepsilon[-3]})\right )
$$ 
where $U''$ is the corresponding open subset of $\mathbb{P}\setminus P$ and $\Phi$ is a specific family of supports, see section 3. Moreover this transform is an isomorphism for suitable open subsets $U$ of $M$. \\

Here is a sketch of the paper. In section 2, we give the relevant materials from sheaf theory that we need in the split Penrose transform. 
In section 3, we review the well-known complex Penrose transform. In sections 4 and 5,  the split Penrose transform is introduced and discussed. In section 6, we finish  the  split Penrose transform and give some examples. Section 7 is devoted to explaining the relationship between our version of the split Penrose transform and the one introduced  in \cite{Ea}. In section 8, we discuss the split Penrose transform in split instanton backgrounds. It is well-known that the Penrose transform has applications in Representation Theory, see \cite{BE}.   
Finally, in section 9, we discuss  the possible applications of the split Penrose transform to Representation Theory (especially representations of $SL(4,\mathbb{R})$).\\
 
\noindent\textbf{Acknowledgment:} I would like to thank Professor I. Frenkel for introducing the subject to me. The author is very grateful to him for support, very useful discussions and his comments on the earlier version of the paper.  I would also like to thank \mbox{Professor} M. Kapranov and \mbox{Professor} G. Zuckerman for helpful and informative conversations.

\end{section}

                      \begin{section}{Preliminaries}
First we review some Cohomology Theory that we will need, see \cite{Br} for the details.\\
Suppose that $X$ is a (Hausdorff) topological space and $\mathcal{F}$ is a sheaf of abelian groups on $X$.  
If $Y$ is a locally closed subset of $X$ (e.g. open or closed), then 
we set $\mathcal{F}_Y$ to be the extension of $\mathcal{F}|_{Y}$ to $X$ by zero. We recall that  the  stalk of   $\mathcal{F}_Y$ at $x$ is $\mathcal{F}_x$  or zero depending on whether $x\in Y$ or $x\notin Y$ respectively.\\
A family of supports on $X$ is a family $\Phi$ of closed subsets of $X$ such that\\
\begin{itemize}
\item a closed subset of a member of $\Phi$ is a member of $\Phi$
\item $\Phi$ is closed under finite unions

\end{itemize}
$\Phi$ is called \textbf{paracompactifying} family of supports if in addition
\begin{itemize}
\item each element of $\Phi$ is paracompact
\item each element of $\Phi$ has a (closed) neighborhood which is in $\Phi$
\end{itemize}
If $S$ is a subspace of $X$ and $\Phi$ is a family of supports on $X$, then it is easy to see that 
$$\Phi|S:=\{A\in\Phi | A\subset S\}$$
is a family of supports on $S$ (and X).\\
Suppose that $\Phi$ is a family of supports on $X$. The set of global sections of $\mathcal{F}$ whose supports are in $\Phi$ is denoted by $\Gamma_{\Phi}(X,\mathcal{F})$. It is well-known that 
the functor $\mathcal{F}\to \Gamma_{\Phi}(X,\mathcal{F})$ is left exact and its $p$-th right derived functor  is denoted by $H^p_{\Phi}(X,\mathcal{F})$. If $S\subset X$ and $\Phi$ is a family of supports on $S$ then we denote 
 $H^n_{\Phi}(S,\mathcal{F}|_S)$ by $H^n_{\Phi}(S,\mathcal{F})$.  \\
Now suppose that $Y$ is a closed subset of $X$. Then it is well-known that one has the following canonical exact sequence of sheaves on $X$
$$0\to\mathcal{F}_{U}\to\mathcal{F}\to\mathcal{F}_Y\to 0$$
where $U=X\setminus Y$.
This exact sequence gives rise to the following exact sequence
$$0\to H^0_{\Phi}(X,\mathcal{F}_U)\to H^{0}_{\Phi}(X,\mathcal{F})\to H^0_{\Phi}(X,\mathcal{F}_Y)$$
$$\quad \qquad\; \to H^1_{\Phi}(X,\mathcal{F}_U)\to H^{1}_{\Phi}(X,\mathcal{F})\to H^1_{\Phi}(X,\mathcal{F}_Y)\to \cdots$$
for any family of supports $\Phi$.
The following lemma characterizes various terms of this exact sequence. 
\begin{lemma}
(a) The natural maps 
$$H^*_{\Phi}(X,\mathcal{F}_{Y})\to H^{*}_{\Phi|Y}(Y,\mathcal{F})$$
are isomorphisms.\\
(b) Suppose that $\Phi$ is a paracompactifying family of supports. Then, there are natural isomorphisms \\
$$H^*_{\Phi}(X,\mathcal{F}_U)\overset{\cong}{\to} H^{*}_{\Phi |U}(U,\mathcal{F})$$
 
\end{lemma}
\begin{proof}
These are standard facts in sheaf theory, see \cite{Br}.
\end{proof}
This lemma has the following corollary
\begin{corollary}\label{LE}
Suppose that  $Y$ is a closed subset of $X$ and $\Phi$ is a paracompactifying family of supports on $X$. Then we have the following canonical exact sequence
of abelian groups
$$0\to H^0_{\Phi |U}(U,\mathcal{F})\to H^{0}_{\Phi}(X,\mathcal{F})\to H^0_{\Phi| Y}(Y,\mathcal{F})$$
$$\qquad \quad\to H^1_{\Phi |U}(U,\mathcal{F})\to H^{1}_{\Phi }(X,\mathcal{F})\to H^1_{\Phi|Y}(Y,\mathcal{F})\to \cdots$$
for any sheaf of abelian groups $\mathcal{F}$ where $U=X\setminus Y$.
\end{corollary}
We use the following conventions when $\Phi$ is one of the following special families of supports. If $\Phi$ is the set of all closed subsets of $X$, then we drop the index $\Phi$. If $\Phi$ is the family of all compact subsets of $X$, we use $c$ instead of $\Phi$. With this terminology, we note that if $X$ is compact and $\Phi$ is the set of all closed subsets of $X$ then 
the exact sequence in corollary  \ref{LE} becomes
\begin{equation}
0\to H^0_c(U,\mathcal{F})\to H^{0}(X,\mathcal{F})\to H^0(Y,\mathcal{F})
\end{equation}
$$\qquad \quad\to H^1_c(U,\mathcal{F})\to H^{1}(X,\mathcal{F})\to H^1(Y,\mathcal{F})\to \cdots$$
 

\noindent \textbf{Twisting Sheaf}\\
Suppose that $X$ is a smooth (or complex) manifold and $\pi:\widetilde{X}\to X$ is a 2-covering of $X$ with $\sigma:\widetilde{X}\to\widetilde{X}$ the involution satisfying $\pi\circ\sigma=\pi$. Then there is  a canonical locally constant sheaf $\widetilde{\mathbb{C}}$ on $X$ such that for any open subset $U\subset X$, $\widetilde{\mathbb{C}}(U)$ is just the set of complex-valued locally constant functions on $\widetilde{U}:=\pi^{-1}(U)$ such that there are odd with respect to $\sigma$.  We call this sheaf the twisting sheaf. In general, if $\mathcal{F}$ is a sheaf of abelian groups on $X$,   then the twisted sheaf of $\mathcal{F}$, denoted by $\widetilde{\mathcal{F}}$, is defined to be the sheaf which consists of sections of $\pi^{-1}\mathcal{F}$  odd with respect to $\sigma$. 
In the same way, if $E$ is a smooth (or holomorphic) vector bundle over $X$ then we can define the twisted vector bundle $\widetilde{E}$.  \\
It is easy to see that, for any sheaf $\mathcal{F}$ on $X$, we have a natural isomorphism
$$\mathcal{F}\oplus\widetilde{\mathcal{F}}\cong\pi_*\pi^{-1}\mathcal{F} $$
(this is just the analog of the decomposition of functions into the sum of even and odd functions). Here $\pi^{-1}\mathcal{F}$ is the inverse image sheaf of $\mathcal{F}$ and $\pi_*$ stands for the direct image functor.
 It is also easy to see that there are natural isomorphisms 
$$ H^n(X,\pi_*\pi^{-1}\mathcal{F})\cong H^n(\widetilde{X},\pi^{-1}\mathcal{F})$$ for any $n$. Therefore we have natural isomorphisms 
\begin{equation}\label{twistedcohomology}
H^n(X,\mathcal{F})\oplus H^n(X,\widetilde{\mathcal{F}})\overset{\cong}{\to}  H^n(\widetilde{X},\pi^{-1}\mathcal{F}) 
\end{equation}
for any $n$.   This is the generalization of the simple fact that every function on $\widetilde{X}$ can be uniquely written as a sum of odd and even functions with respect to $\sigma$. \\
\noindent \textbf{Conventions}: We use $\varepsilon$, $\mathcal{O}$ and $\omega$ for smooth, holomorphic and real analytic objects. For example, if $V$ is a real analytic vector bundle on a real analytic manifold $X$, we use 
$\Gamma(X,\varepsilon_V)$ to denote the set of global smooth sections of $V$ and $\Gamma^{\omega}(X,\varepsilon_V)$ to denote the set of global real analytic sections of $V$.

\end{section}

                                                                   \begin{section}{Review of the Complex Penrose Transform}
                                        First, we recall the definitions of the basic geometric spaces on which the complex Penrose transform operates. We refer to
  \cite{EPW} and \cite{WW} for more details on the geometry of twistors. \\
The vector space of twistors $\mathbb{T}$ is by definition a four dimensional complex vector space endowed with an Hermitian form $\Phi$ of type $(+\quad+\quad-\quad-)$.  We have the fundamental twistor diagram 

\begin{displaymath}
\xymatrix{ & \mathbb{F} \ar[dl]_{\mu} \ar[dr]^{\nu} & \\
                 \mathbb{P}&  & \mathbb{M}}
 \end{displaymath}
where $\mathbb{P}$ is the space of complex lines in $\mathbb{T}$, $\mathbb{M}$ is the Grassmannian of complex 2-planes in $\mathbb{T}$, $\mathbb{F}$ is the space of pairs of nested 1- and 2-dimensional subspaces of $\mathbb{T}$, and  where $\mu$ and $\nu$ are the natural holomorphic maps.  Both maps $\mu$ and $\nu$ are fiber bundle maps where the fibers of $\mu$ are isomorphic to $\mathbb{CP}^2$ and the fibers of $\nu$ are isomorphic to $\mathbb{CP}^1$.  Given any open subset $U$ of 
$\mathbb{M}$ we set
$$U':=\nu^{-1}(U), \quad U'':=\mu(U')$$
which are open subspaces of $\mathbb{F}$ and $\mathbb{P}$ respectively.  Suppose that $V$ is a holomorphic vector bundle on $\mathbb{P}$ and $\mathcal{O}_V$ is the sheaf of holomorphic sections of $V$. One can summarize the complex Penrose transform in three steps as follow: (see \cite{EPW})\\
(1) \textit{Pull-back step:}\\ 
There are natural maps 
$$\mu^*:H^n(U'',\mathcal{O}_V)\to H^n(U',\mu^{-1}\mathcal{O}_V)$$
We recall that the map $\mu:U'\to U''$  is called \textbf{elementary} if its fibers are connected and have vanishing first Betti number. If  
the map $\mu:U'\to U''$ is elementary then 
$$\mu^*:H^0(U'',\mathcal{O}_V)\to H^0(U',\mu^{-1}\mathcal{O}_V)$$
and
$$\mu^*:H^1(U'',\mathcal{O}_V)\to H^1(U',\mu^{-1}\mathcal{O}_V)$$
are bijections.  \\
 (2) \textit{Middle step:}\\ 
There is an exact sequence of sheaves on $\mathbb{F}$
$$0\longrightarrow\mu^{-1}\mathcal{O}_V\longrightarrow\mathcal{O}_{\mu^*V}\overset{d_\mu}{\longrightarrow}\Omega^1_\mu(V)\overset{d_\mu}{\longrightarrow}\Omega^2_\mu(V)\longrightarrow 0$$
where $\Omega^n_\mu(V)=\mathcal{O}_{\mu^*V}\otimes\Omega^n_{\mu}$ ($\Omega^n_{\mu}$ is the sheaf of holomorphic relative $n$-forms on $\mathbb{F}$ with respect to the fibration $\mu$), $d_{\mu}$ is the induced exterior derivative on forms and $\mu^*V$ is the pull-back vector bundle.  This exact sequence gives rise to a spectral sequence. More precisely there is a spectral sequence
$$E^{p,q}_1=H^q(U',\Omega^p_\mu(V))\Longrightarrow H^{p+q}(U',\mu^{-1}\mathcal{O}_{V})$$
where the differentials $d_1:E_1^{p,q}\to E^{p+1,q}_1$ are induced by the relative exterior derivative
$d_\mu:\Omega^p_\mu(V)\to\Omega^{p+1}_\mu(V)$.\\
 (3) \textit{Push-forward step:} \\
The  Leray spectral sequence of $\nu:U'\to U$ relates $H^q(U',\Omega^p_\mu(V))$ with cohomology groups on $U$ with coefficients in the direct image sheaves $v^q_{*}\Omega^p_\mu(V)$. \\ \\

For appropriate holomorphic vector bundles $V$ on $\mathbb{P}$, these steps give a map from $H^1(U'',\mathcal{O}_V)$ to the  (co)kernel of differential operators between vector bundles on $U$.  We introduce the vector bundles and differential operators on $\mathbb{M}$ which appear in the Penrose transform, see \cite{EPW} (we make no distinction between primed and unprimed spinor bundles).  We  denote the universal vector bundle of $\mathbb{M}$ by $H$.  We use the following notations
$$H_n:= \text{the}\; n\text{-th  symmetric product of} \; H$$
$$\mathcal{O}[-1]:= H\wedge H$$
$$\mathcal{O}[1]:=\text{the dual of}\;\mathcal{O}[-1]$$
  $$\mathcal{O}[k]:=\otimes^{k}\mathcal{O}[1] \;\text{for} \;k\in\mathbb{Z} 
$$
Finally, for any vector bundle $E$ on $\mathbb{M}$, we denote $E\otimes\mathcal{O}[k]$ by $E[k]$. 
For any $n\geq 1$, there are  first order linear differential operators (see \cite{EPW})
$$D_n:\Gamma(\mathbb{M},\mathcal{O}_{H_n[-1]})\to \Gamma(\mathbb{M},\mathcal{O}_{H\otimes H_{n-1}[-2]})$$
We denote the kernel of $D_n$ on an open subset $U$ of $M$  by $\mathcal{Z}_n(U)$ which is the set of "holomorphic massless fields on $U$ of helicity $n/2$". Finally we have the Laplacian operator (wave operator)
$$\square:\Gamma(\mathbb{M},\mathcal{O}[-1])\to\Gamma(\mathbb{M},\mathcal{O}[-3])$$
We also denote the kernel of $D_0$ on an open subset $U$ of $M$  by $\mathcal{Z}_0(U)$.\\
Using these vector bundles and differential operators, we can give the Penrose transform applied to $\mathcal{O}(-n-2)$ ($n \geq 0$). More precisely, we have a  transform 
\begin{equation}\label{PT}
H^1(U'',\mathcal{O}(-n-2))\overset{\mathcal{P}}{\longrightarrow}\mathcal{Z}_{n}(U)
\end{equation}
Moreover,  if 
the map $\mu:U'\to U''$ is elementary then this transformation is  a bijection.   \\ \\
\begin{remark}\label{VM}
We note that $M=Gr(2,\mathbb{R}^4)$ is the real counterpart of $\mathbb{M}=Gr(2,\mathbb{C}^4)$. Therefore all the canonical vector bundles on $\mathbb{M}$ and natural differential operators between them as above can be defined on $M$ completely similarly (which are, in fact, the restrictions to $M$) .  We use the same conventions to denote these objects on $M$. For example we denote the counterpart of $\mathcal{O}[-1]$ on $M$ by $\varepsilon[-1]$ and  the counterpart of the Laplacian is the ultrahyperbolic operator 
 $$\square_{2,2}:\Gamma(
M,\widetilde{\varepsilon[-1]})\to\Gamma(M,\widetilde{\varepsilon[-3]})$$
 see \cite{Ea,BaE}.

\end{remark}
 Following the complex Penrose transform we would like to present the split Penrose transform in three steps. As we will see, in our approach to the split Penrose transform, the second and third steps are essentially the steps in the complex Penrose transform combined with a direct limit process.  In the following sections we introduce the split Penrose transform.

 \end{section}
\begin{section}{ Split Twistor Diagram}               
We would like to obtain a split version of the Penrose transform. The complex Penrose transform is 
an $SL(4,\mathbb{C})$-equivariant transform. Different real forms of $SL(4,\mathbb{C})$ lead to different versions of the Penrose transform. Here we consider the real form $SL(4,\mathbb{R})$. This amounts to 
consider $M$, the Grassmannian of 2-planes in $\mathbb{R}^4$, as a totally real submanifold of $\mathbb{M}$. By the split Penrose transform we mean a transform which
 identifies cohomological data on $\mathbb{P}$ (or appropriate open subsets of it) with the solutions of relevant differential equations on $M$.\\
First we introduce the appropriate geometric setting for the transformation in the split signature. 
Suppose that $T$ is a real vector subspace of $\mathbb{T}$ of real dimension 4.  Then we have the  following "real" analog of the twistor diagram (see \cite{GS})
   
\begin{displaymath}   
\xymatrix{ & F \ar[dl]_{\mu_0} \ar[dr]^{\nu_0} & \\
                 P&  & M}
 \end{displaymath}  
 where $P$  is the space of real lines in $T$, $M$ is the Grassmannian of real 2-planes in $T$, $F$ is the space of pairs of nested 1- and 2-dimensional subspaces of $T$ and  where $\mu_0$ and $\nu_0$ are the natural projection maps. Both maps $\mu_0$ and $\nu_0$ are fiber bundle maps where the fibers of $\mu_0$ are isomorphic to $\mathbb{RP}^2$ and the fibers of $\nu_0$ are isomorphic to $\mathbb{RP}^1$. In fact $\mathbb{P}$, $\mathbb{F}$ and $\mathbb{M}$ are the complexifications of $P$, $F$ and $M$  respectively. \\
 Both manifolds $M$ and $F$ have double covers $\widetilde{M}$ and $\widetilde{F}$. More precisely,
 $\widetilde{M}$ is the Grassmannian of \emph{oriented} 2-planes in $T$ and $\widetilde{F}$ is the space of pairs of nested 1-dimensional  and  \emph{oriented} 2-dimensional subspaces of $T$.  This gives us the following diagram  
 \begin{displaymath}   
\xymatrix{ & \widetilde{F} \ar[dl]_{\widetilde{\mu_0}} \ar[dr]^{\widetilde{\nu_0}} & \\
                 P&  & \widetilde{M}}
 \end{displaymath}  
We can consider the corresponding twisting sheaves on $M$ and $F$. Moreover, we note that the pull back of the twisting sheaf on $M$ via $\nu_0$ is just the twisting sheaf on $F$.   We see that this double coverings and consequently the corresponding twisting sheaves have an extension to appropriate open subsets of $\mathbb{M}$ and $\mathbb{F}$. More precisely, consider a non-degenerate symmetric bilinear form $(\;,\:)$ on $\mathbb{T}$ coming from an inner product on $T$. Let $\mathbb{M}_r$ be the set of planes $p\in \mathbb{M}$ which the restriction of $(\;,\:)$ on them is non-degenerate. It is easy to see that $\mathbb{M}_r$ is an open subset of $\mathbb{M}$ which contains $M$. Moreover there is a natural double cover $\widetilde{\mathbb{M}_r}$ of  $\mathbb{M}_r$ coming from $(\:,\:)$ whose restriction to $M$ is just its natural double cover $\widetilde{M}$. In the same way, $\mathbb{F}_r:=\nu^{-1}(\mathbb{M}_r)$ has a natural double cover $\widetilde{\mathbb{F}_r}$. Moreover $\nu$ extends to a canonical map $\nu: \widetilde{\mathbb{F}_r}\to\widetilde{\mathbb{M}_r}$. \\
Following the complex case, for any open subset $U$ of $M$ we define the following open sets
$$U'_R:=\nu_0^{-1}(U),  \quad U''_R:=\mu_0(U'_R) $$
We note that $\widetilde{U'_R}=\widetilde{\nu_0}^{-1}(\widetilde{U})$. \\
 The last space which enters into our picture is the following space 
$$G:=\nu^{-1}(M)$$
In other words, $G$ is the space of pairs of complex lines $L$ and real 2-planes $K$ such that $L\subset K\otimes\mathbb{C}$.
It is easy to see that $F\subset G\subset \mathbb{F}$. By abuse of notation we denote the restriction of $\mu:\mathbb{F}\to \mathbb{P}$ to $G$ by $\mu:G\to \mathbb{P}$ as well.   We have the following simple lemma (see  \cite{Ea} or \cite{Ma})
\begin{lemma}\label{diffeomorphism}
The map $\mu:G\setminus F\to \mathbb{P}\setminus P$ is a diffeomorphism. 
\end{lemma} 
 
We note that $G$ also has a double cover $\widetilde{G}$ which is just the space of pairs of complex lines $L$ and \emph{oriented} real 2-planes $K$ such that $L\subset K\otimes\mathbb{C}$. Therefore there is a twisting sheaf on $G$ coming from this double cover. We note that we have the inclusions $F\subset G\subset \mathbb{F}_r$ and natural inclusions 
$\widetilde{F}\subset \widetilde{G}\subset \widetilde{\mathbb{F}_r}$ and hence the twisting sheaves are compatible. Moreover we have (see \cite{Ea})
\begin{lemma}\label{DT}
 For any sheaf of abelian groups $\mathcal{F}$ on $G$, we have canonical isomorphisms
 $$ \nu^n_{*}\widetilde{\mathcal{F}}\cong\widetilde{\nu^n_{*}\mathcal{F}}$$
 where $\nu:G\to M$ is the restriction of $\nu:\mathbb{F}\to\mathbb{M}$ and $\nu^n_*$ is the $n$-th direct image functor.
 \end{lemma}
 Finally, we define the split twistor diagram to be the following diagram (see \cite{Ea})
 \begin{displaymath}  
\xymatrix{ & G \ar[dl]_{\mu} \ar[dr]^{\nu} & \\
                 \mathbb{P}&  & M}
 \end{displaymath}  
 where $\mu$ and $\nu$ are the obvious maps as always.  Like the complex case, for an open subset $U$ of $M$ we also define
 $$U':=\nu^{-1}(U)$$

 \end{section}

 \begin{section}{Split Penrose Transform: the pull-back step}
In this section we provide the first step of the split Penrose transform. In the first step, given an open subset $U$ of $M$, we would like to start with a first cohomology group on an appropriate open subset of $\mathbb{P}$ and  identify it with a first cohomology group on $U'$. First of all, it is not clear which open set of $\mathbb{P}$ to take. Moreover, it turns out that the usual Cohomology Theory would not suffice and we will need to use Cohomology Theory with supports. As we will see later,
the right open set of $\mathbb{P}$ corresponding to $U$ is 
$$U'':=\mu(U'\setminus F)$$
The family of supports on $U''$ which appears in the split Penrose transform is the following
$$\Phi(U):=\{ A\subset U'' | \quad \mu^{-1}(A) \quad \text{ is a closed subset of} \quad U' \}$$
\begin{remark}\label{fibration}
Note that there is a  canonical map $\pi:\mathbb{P}\setminus P\to M$ sending a complex line in $\mathbb{C}^4$ to the plane generated by its real and imaginary parts in $\mathbb{R}^4$. Then it is easy to see that 
$U''=\pi^{-1}(U')$.
\end{remark}
For the pull-back step we need the following proposition
\begin{proposition}\label{PB}
Suppose that $X$ and $Y$ are real analytic manifolds and $f:X\to Y$ is a surjective smooth map of maximal rank. If $f:X\to Y$ is elementary, i.e. the fibers are connected and have vanishing first Betti number, then the natural maps 
$$f^*:H^n(Y,\omega_V)\to H^n(X,f^{-1}\omega_V)$$
are isomorphisms for $n=0,1$ and any real analytic vector bundle $V$ on $Y$. Here $\omega_V$ denotes the sheaf of  real analytic sections of $V$.

\end{proposition}
\begin{proof}
It is easy to see that, since the fibers of $f$ are connected, for any sheaf $\mathcal{F}$ on $Y$ the natural map 
$$f^*:H^0(Y,\mathcal{F})\to H^0(X,f^{-1}\mathcal{F})$$
is an isomorphism. On the other hand, we note that  these hypotheses imply that $H^1(X,f^{-1}\varepsilon_V)=0$, see \cite{Bu}.  The following exact sequence 
$$0\to\omega_V\to\varepsilon_V\to\frac{\varepsilon_V}{\omega_V}\to 0$$
gives rise to the following exact sequence
$$0\to H^0(Y,\omega_V)\to H^0(Y,\varepsilon_V)\to H^0(Y,\frac{\varepsilon_V}{\omega_V})\to  H^1(Y,\omega_V) \to 0$$
because $H^1(Y,\varepsilon_V)=0$. Similarly, the exact sequence 
$$0\to f^{-1}\omega_V\to f^{-1}\varepsilon_V\to f^{-1}\frac{\varepsilon_V}{\omega_V}\to 0$$
gives the following exact sequence
$$0\to H^0(X,f^{-1}\omega_V)\to H^0(X,f^{-1}\varepsilon_V)\to H^0(X,f^{-1}\frac{\varepsilon_V}{\omega_V})\to  H^1(X,f^{-1}\omega_V) \to 0$$
because $H^1(X,f^{-1}\varepsilon_V)=0$.  Now consider the following commutative diagram
\begin{displaymath}  
\xymatrix{ 0\ar[r] & H^0(Y,\omega_V)\ar[d]^{f^*} \ar[r] & H^0(Y,\varepsilon_V) \ar[d]^{f^*}\ar[r] & H^0(X,f^{-1}\frac{\varepsilon_V}{\omega_V})\ar[d]^{f^*}\ar[r] & H^1(Y,\omega_V)\ar[d]^{f^*}\ar[r] & 0\\
0\ar[r] & H^0(X,f^{-1}\omega_V)\ar[r] & H^0(X,f^{-1}\varepsilon_V)\ar[r] & H^0(X,f^{-1}\frac{\varepsilon_V}{\omega_V})\ar[r] & H^1(X,f^{-1}\omega_V)\ar[r] & 0}
 \end{displaymath} 
The rows are exact and the first three vertical maps are isomorphisms. Therefore the last one is also an isomorphism.
\end{proof}
Now we can give the pull-back step of the split Penrose transform.
\begin{theorem}\label{MT}
Let $V$ be a holomorphic vector bundle on $\mathbb{P}$. Then, for any open subset $U$ of $M$ there is a canonical map (defined up to sign)  
$$H^1_{\Phi(U)}(U'',\mathcal{O}_V)\overset{\mathcal{P}_0}{\to} H^1(U',\widetilde{\mu^{-1}\mathcal{O}_V})$$
Moreover, if the maps $\mu_0:U''_R\to U'_R$ and  $\widetilde{\mu_0}:\widetilde{U''_R}\to U'_R$ are elementary then $\mathcal{P}_0$ is an isomorphism.
\end{theorem}
\begin{proof}
Clearly the natural map
$$\mu^*: H^1_{\Phi(U)}(U'',\mathcal{O}_V)\to H^1_{\mu^{-1}\Phi(U)}(U'\setminus F,\mu^{-1}\mathcal{O}_V)$$ 
is an isomorphism, see lemma \ref{diffeomorphism}.  Here $\mu^{-1}\Phi(U)$ is defined to be the family of sets $\mu^{-1}A$ where $A\in \Phi(U)$. Since the twisting sheaf on $G\setminus F$ is trivial, we have a canonical isomorphism (up to sign)
$$H^1_{\mu^{-1}\Phi(U)}(U'\setminus F,\mu^{-1}\mathcal{O}_V)\to H^1_{\mu^{-1}\Phi(U)}(U'\setminus F,\widetilde{\mu^{-1}\mathcal{O}_V})$$
Using corollary \ref{LE}, we obtain the following exact sequence
\begin{equation}\label{ES}
H^0(U'_R,\widetilde{\mu_0^{-1}\mathcal{O}_V})\to H^1_{\mu^{-1}\Phi(U)}(U'\setminus F,\widetilde{\mu^{-1}\mathcal{O}_V})\to H^1(U',\widetilde{\mu^{-1}\mathcal{O}_V})\to H^1(U'_R,\widetilde{\mu_0^{-1}\mathcal{O}_V})
\end{equation}
where the family is chosen to be the set of all closed subsets of $U'$.
So the middle map  composed with the previous maps provides the desired map
$$H^1_{\Phi(U)}(U'',\mathcal{O}_V)\overset{\mathcal{P}_0}{\to} H^1(U',\widetilde{\mu^{-1}\mathcal{O}_V})$$
Now suppose that the maps $\mu_0:U''_R\to U'_R$ and  $\widetilde{\mu_0}:\widetilde{U''_R}\to U'_R$ are elementary. Then from  proposition \ref{PB}, both maps 
$$H^0(U''_R,\mathcal{O}_V){\to}H^0(U'_R,\mu_0^{-1}\mathcal{O}_V) $$ 
$$H^1(U''_R,\mathcal{O}_V){\to}H^1(\widetilde{U'_R},\widetilde{\mu_0}^{-1}\mathcal{O}_V) $$ 
are isomorphisms.  This implies that the natural maps 
$$H^n(U'_R,\mu_0^{-1}\mathcal{O}_V) \to H^n(\widetilde{U'_R},\widetilde{\mu_0}^{-1}\mathcal{O}_V) $$
are isomorphisms for $n=0,1$. By isomorphism \ref{twistedcohomology}, we see that 
$$H^n(U'_R,\widetilde{\mu_0^{-1}\mathcal{O}_V})=0$$ 
for $n=0,1$. This implies that 
the middle map in the exact sequence \ref{ES} is an isomorphism. Hence the map 
$$H^1_{\Phi(U)}(U'',\mathcal{O}_V)\overset{\mathcal{P}_0}{\to} H^1(U',\widetilde{\mu^{-1}\mathcal{O}_V})$$
is an isomorphism.

\end{proof}
The importance of this theorem is that even though we start with cohomology groups on  open subsets other than the ones that the usual Penrose transform suggests and we work with Cohomology Theory with supports, after passing to $G$ we obtain the cohomology groups  which can be handled by the complex Penrose transform as we will see in the next section.  
\end{section}

\begin{section}{Split Penrose Transform: second and third steps and examples} 
 As we saw in the last section there is a canonical (up to sign) map
 $$H^1_{\Phi(U)}(U'',\mathcal{O}_V)\to H^1(U',\widetilde{\mu^{-1}\mathcal{O}_V})$$
 for any open subset $U$  of $M$. Now we want to transform this cohomological data on $U'$ down to $U$. This can be done by a direct limit process. More precisely, suppose that $\{\mathbb{U}_i\}$ is a fundamental system of neighborhoods of $U$ in $\mathbb{M}_r$.  The fibers of $\nu:\mathbb{F}\to\mathbb{M}$ are compact, so $\{\mathbb{U}_i'=\nu^{-1}(\mathbb{U}_i)\}$ forms a fundamental system of neighborhoods for $U'$. Therefore  the canonical map
$$ \underset{i}{\underrightarrow{\lim}}\: H^1(\mathbb{U}_i',\widetilde{\mu^{-1}\mathcal{O}_V})\to H^1(U',\widetilde{\mu^{-1}\mathcal{O}_V})$$
 is an isomorphism, see \cite{Br}. Now we can invoke the complex Penrose transform to transform the terms in this direct limit down to $\mathbb{M}$ (the twisting sheaf introduces no difficulties  by lemma \ref{DT}). This transformation clearly commutes with the direct limit and hence we obtain data on $U$.  In other words, the second and third steps of the split Penrose transform are basically those in the complex Penrose transform. Now,  we give some examples to clarify the split Penrose transform.\\
The first important class of examples is obtained by taking $V$ to be a line bundle. Suppose that $\mathcal{O}_V=\mathcal{O}(-n-2)$ where $n\geq 0$.  Then the split Penrose transform gives the following theorem
\begin{theorem}
 Let  $U$ be an open subset of $M$. There is a  map
$$H^1_{\Phi(U)}(U'',\mathcal{O}(-n-2))\overset{P}{\to}\widetilde{\mathcal{Z}_n(U)}^\omega$$
where  $\widetilde{\mathcal{Z}_{n}(U)}^\omega$ is the set of real analytic twisted massless fields on $U$ of helicity $n/2$, see remark \ref{VM}. Moreover if  the maps $\mu_0:U''_R\to U'_R$ and  $\widetilde{\mu_0}:\widetilde{U''_R}\to U'_R$ are elementary then the above map is an isomorphism. 
 \end{theorem}
\begin{proof}
By theorem \ref{MT} we have a map
$$H^1_{\Phi(U)}(U'',\mathcal{O}(-n-2))\overset{\mathcal{P}_0}{\to} H^1(U',\widetilde{\mu^{-1}\mathcal{O}(-n-2)})$$
By the above discussion we have an isomorphism
$$ \underset{i}{\underrightarrow{\lim}}\: H^1(\mathbb{U}_i',\widetilde{\mu^{-1}\mathcal{O}(-n-2)})\to H^1(U',\widetilde{\mu^{-1}\mathcal{O}(-n-2)})$$
Now the complex Penrose transform (steps 2 and 3) gives an isomorphism
$$H^1(\mathbb{U}_i',\widetilde{\mu^{-1}\mathcal{O}(-n-2)})\to \widetilde{\mathcal{Z}_n(\nu(\mathbb{U}_i))}$$
see transform \ref{PT} and lemma \ref{DT}.  Therefore we obtain a map 
$$H^1_{\Phi(U)}(U'',\mathcal{O}(-n-2))\to \underset{i}{\underrightarrow{\lim}}\:\widetilde{\mathcal{Z}_n(\nu(\mathbb{U}_i))}$$
But it is clear that 
$$\widetilde{\mathcal{Z}_n(U)}^\omega\cong\underset{i}{\underrightarrow{\lim}}\:\widetilde{\mathcal{Z}_n(\nu(\mathbb{U}_i))}$$
Therefore we obtain the desired transform. If the maps $\mu_0:U''_R\to U'_R$ and  $\widetilde{\mu_0}:\widetilde{U''_R}\to U'_R$ are elementary, then by theorem \ref{MT}, the transform $\mathcal{P}$ is an isomorphism.

\end{proof}
Note that for $n=0$ we obtain a map 
 $$H^1_{\Phi(U)}(U'',\mathcal{O}(-2))\overset{\mathcal{P}}{\to}\ker\left (\square_{2,2}:\Gamma^{\omega}(
U,\widetilde{\varepsilon[-1]})\to
\Gamma^{\omega}(U,\widetilde{\varepsilon[-3]})\right )
$$ 
This is a local version of the X-ray transform. The global picture is also interesting. More precisely, if we take $U=M$, then it is easy to see that $U'=G$, $U''_R=P$ and $U''=\mathbb{P}\setminus P$. Moreover the maps $\mu_0:U''_R\to U'_R$ and  $\widetilde{\mu_0}:\widetilde{U''_R}\to U'_R$ are elementary. Therefore, the first step of the split Penrose transform provides an isomorphism. Finally, it is easy to see that $\Phi(M)$ is just the set of all compact subsets of $\mathbb{P}\setminus P$. Therefore we obtain the following isomorphism
$$H^1_c(\mathbb{P}\setminus P,\mathcal{O}(-2))\overset{\mathcal{P}}{\to}\ker\left (\square_{2,2}:\Gamma^{\omega}(
M,\widetilde{\varepsilon[-1]})\to
\Gamma^{\omega}(M,\widetilde{\varepsilon[-3]})\right )
$$
This is the cohomological interpretation of the X-ray transform as we will see in the next section. \\

\end{section}

\begin{section}{Comparison Between Two Versions of the Split Penrose Transform}

The first version of the split Penrose transform was appeared in \cite{Ea}. In this section we explain the relationship between the two versions of the split Penrose transforms.   We must point out that our version of the split Penrose transform deals with  real analytic objects whereas the version in \cite{Ea} deals with smooth ($C^{\infty}$) objects.\\
 
First we recall the split Penrose transform as in \cite{Ea} which we call the "\textit{smooth}" Penrose transform to distinguish it from our version. In the smooth Penrose transform one starts with smooth data on $\mathbb{RP}^3$ (namely smooth sections of vector bundles) and identify that with the (co)kernel of certain differential operators on $M$.
The diagram for the smooth split Penrose transform is also
\begin{displaymath}  
\xymatrix{ & G \ar[dl]_{\mu} \ar[dr]^{\nu} & \\
                 \mathbb{P}&  & M}
 \end{displaymath} 
 The main step in the smooth Penrose transform is to interpret  smooth data on $\mathbb{RP}^3$ as some data on $G$. Even though $G$ is not a complex manifold, it is very close to be one (see lemma \ref{diffeomorphism}). More precisely, there is an "involutive" structure on $G$. This involutive structure can be used to define the involutive cohomology, for the details see \cite{Ea}. The key theorem in the smooth Penrose transform is 
 \begin{theorem}\label{Ea}
 For any holomorphic vector bundle $V$ on $\mathbb{P}$ there is an exact sequence
 $$0\to\Gamma(\mathbb{P},\mathcal{O}_V)\to\Gamma(P,V)\to H^1_{in}(G,\widetilde{\mu^*V})\to H^1(\mathbb{P},\mathcal{O}_V)\to 0$$
 where $\Gamma(\mathbb{P},\mathcal{O}_V)$ is the set of global holomorphic sections of $V$, $\Gamma(P,V)$ is the set of global smooth sections of $V|_P$ and $H^1_{in}(G,\widetilde{\mu^*V})$ is the first "involutive cohomology" group associated to $\widetilde{\mu^*V}$. 
 \end{theorem}
 Therefore, up to finite dimensional spaces,  $\Gamma(P,V)\to H^1_{in}(G,\widetilde{\mu^*V})$ is an isomorphism. Now the next step in the smooth Penrose transform is to interpret $H^1_{in}(G,\widetilde{\mu^*V})$ as smooth data on $M$. 
Due to the introduction of the involutive cohomology, the other steps of the smooth Penrose transform are more involved. Nevertheless, similar to the complex Penrose transform, one can identify 
$H^1_{in}(G,\widetilde{\mu^*V})$ as smooth data on $M$ (note that the fibers of $\mu:G\to M$ are copies of $\mathbb{CP}^1$), for the full description of the smooth Penrose transform see \cite{Ea,BaE}  .\\  
As an example of the smooth Penrose transform, one has the so-called X-ray (or sometimes called Radon) transform
$$\Gamma(P
,\mathcal{\varepsilon}(-2))\overset{\mathcal{R}}{\to}\ker\left (\square_{2,2}:\Gamma(
M,\widetilde{\varepsilon[-1]})\to
\Gamma(M,\widetilde{\varepsilon[-3]})\right )
$$
which is a bijection, see \cite{He} for a full discussion on the Radon transform. We recall that  $\varepsilon(-1)$ is the universal line bundle on $P=\mathbb{RP}^3$ and $\varepsilon(-2)=\varepsilon(-1)\otimes\varepsilon(-1)$. \\
As one might guess, the split Penrose transform must be hidden in the smooth Penrose transform.  First of all we have the following theorem,
\begin{theorem}\label{relation}
For any holomorphic vector bundle $V$ on $\mathbb{P}$, there is a natural one-to-one map 
$j:H^1_{c}(\mathbb{P}\setminus P,\mathcal{O}_V)\to H^1_{in}(G,\widetilde{\mu^*V})$ for which the following diagram is commutative
 \begin{displaymath}  
\xymatrix{ 0\ar[r] & \Gamma(\mathbb{P},\mathcal{O}_V)\ar@{=}[d] \ar[r] & \Gamma^{\omega}(P,V) \ar[d]^{i}\ar[r] & H^1_{c}(\mathbb{P}\setminus P,\mathcal{O}_V)\ar[d]^{j}\ar[r] & H^1(\mathbb{P},\mathcal{O}_V)\ar@{=}[d]\ar[r] & 0\\
0\ar[r] & \Gamma(\mathbb{P},\mathcal{O}_V)\ar[r] & \Gamma(P,V) \ar[r] & H^1_{in}(G,\widetilde{\mu^*V})\ar[r] & H^1(\mathbb{P},\mathcal{O}_V)\ar[r] & 0}
 \end{displaymath} 
 where $i: \Gamma^{\omega}(P,V)\to \Gamma(P,V)$  is just the inclusion map.

\end{theorem}
\begin{proof} We recall the construction of the exact sequence in theorem \ref{Ea}. Consider the Dolbeault complex associated to $V$, 
$$0\to\Gamma(\mathbb{P},V)\overset{\bar{\partial}}{\to}\Gamma(\mathbb{P},\varepsilon^{0,1}_V)\overset{\bar{\partial}}{\to}\Gamma(\mathbb{P},\varepsilon^{0,2}_V)\overset{\bar{\partial}}{\to}\Gamma(\mathbb{P},\varepsilon^{0,3}_V)\to 0
$$
where $\Gamma(\mathbb{P},\varepsilon^{0,i}_V)$ is the set of smooth $V$-valued $(0,i)$-forms on $\mathbb{P}$. We denote this complex by $\Gamma(V)$. Then we can consider the sub-complex $\Gamma_{\infty}(V)$ of $\Gamma(V)$ consisting of  sections which are zero along $P$ to infinite order, see \cite{Ea}.  Then one considers the following exact sequence of complexes
$$0\to\Gamma_{\infty}(V)\to\Gamma(V)\to\frac{\Gamma(V)}{\Gamma_{\infty}(V)}\to 0$$
The corresponding exact sequence of the cohomology groups yields  the exact sequence in theorem \ref{Ea}
 $$0\to\Gamma(\mathbb{P},\mathcal{O}_V)\to\Gamma(P,V)\to H^1_{in}(G,\widetilde{\mu^*V})\to H^1(\mathbb{P},\mathcal{O}_V)\to 0$$
We can go further and consider a sub-complex of   $\Gamma_{\infty}(V)$. Set  $\Gamma_{0}(V)$
to be the sub-complex of $\Gamma_{\infty}(V)$ consisting of sections which are zero in a neighborhood of $P$.  Since any holomorphic section of $V$ which is zero in a neighborhood of $P$ is identically zero, the corresponding exact sequence of cohomology groups of the  exact sequence 
$$0\to\Gamma_{0}(V)\to\Gamma(V)\to\frac{\Gamma(V)}{\Gamma_{0}(V)}\to 0$$
yields
$$0\to \Gamma(\mathbb{P},\mathcal{O}_V)\to \Gamma^{\omega}(P,V) \to H^1_{c}(\mathbb{P}\setminus P,\mathcal{O}_V)\to H^1(\mathbb{P},\mathcal{O}_V)\to H^1(P,\mathcal{O}_V)$$

We claim that $H^1(P,\mathcal{O}_V)=0$. In fact it is known that any compact totally real submanifold $Y$ of a complex manifold $X$ is a holomorphic set (i.e. there are open Stein submanifolds $Y\subset \dots\subset S_2\subset S_1$ of $X$ which form a fundamental system  of neighborhoods  for $Y$), see \cite{We1}. This, in particular, implies that for any coherent analytic sheaf $\mathcal{F}$ on $X$ we have $H^n(Y,\mathcal{F})=0$ if $n>0$. Clearly $P$ is a compact totally real submanifold of $\mathbb{P}$ and $\mathcal{O}_V$ is a coherent analytic sheaf on $\mathbb{P}$ and hence $H^1(P,\mathcal{O}_V)=0$. Therefore we obtain the following exact sequence
$$0\to \Gamma(\mathbb{P},\mathcal{O}_V)\to \Gamma^{\omega}(P,V) \to H^1_{c}(\mathbb{P}\setminus P,\mathcal{O}_V)\to H^1(\mathbb{P},\mathcal{O}_V)\to 0$$
The inclusion map $\Gamma_{0}(V)\to\Gamma_{\infty}(V)$ provides the maps between the above  exact sequences. Moreover a simple diagram chasing shows that $j$ is injective.

\end{proof}
Using map $j$, we can consider $H^1_{c}(\mathbb{P}\setminus P,\mathcal{O}_V)$ as a subspace of $H^1_{in}(G,\widetilde{\mu^*V})$.  
Therefore, when applying the smooth Penrose transform, one can   keep track of $H^1_{c}(\mathbb{P}\setminus P,\mathcal{O}_V)$.   In particular for $\mathcal{O}(-2)$ one can see the compatibility between the split Penrose transform and the X-ray transform. More precisely, we have the following cohomological description of the X-ray transform
\begin{proposition}
 There is a natural isomorphism
$$H^1_c(\mathbb{P}\setminus P,\mathcal{O}(-2))\cong
\Gamma^{\omega}(P,\varepsilon(-2))$$
and hence an injection
$$H^1_c(\mathbb{P}\setminus P,\mathcal{O}(-2))\to
\Gamma(P,\varepsilon(-2))$$
Moreover the following diagram is commutative 
\begin{displaymath}
\xymatrix{ H^1_c(\mathbb{P}\setminus P,\mathcal{O}(-2)) \ar[d]_{\mathcal{P}}\ar[r] &\Gamma(P,\varepsilon(-2))\ar[d]^{\mathcal{R}} \\
        \ker\square_{2,2}^{\omega} \ar[r]& \ker\square_{2,2}
        }
 \end{displaymath}
where $\mathcal{R}$ is the X-ray transform.

\end{proposition}

\end{section}

\begin{section}{Split Penrose Transform in Split Instanton Backgrounds}
In this section we explain the split Penrose transform in the presence of split instantons, see \cite{We1,Hi} for the Euclidean case.
First we need to explain what we mean by split instantons.  
Here is a short review of the  SDYM equations in the split signature. The space of 2-forms on $M$ has a natural decomposition into two subspaces of self-dual and anti-self-dual 2-forms, see \cite{GS}. For a given Lie group $G$, a $G$-SDYM field on $M$ is a vector bundle $V$ with structure group $G$ and a connection $\nabla$ on $V$ compatible with $G$ whose curvature is self-dual, see \cite{Ar,Ma,MW}.  By a split instanton, we mean a $U(n)$-SDYM field on $M$.  A split instanton $(V,\nabla)$ is called real analytic if both $V$ and $\nabla$ are real analytic.\\
It is known that there is  a canonical holomorphic vector bundle $E$ on $\mathbb{P}=\mathbb{CP}^3$ associated to any split instanton $(V,\nabla)$ on $M$, see \cite{Ma}.
It is easy to describe $E$ on $\mathbb{P}\setminus P$. Consider the fibration $\pi:\mathbb{P}\setminus P\to M$ as in remark \ref{fibration}. Then it is easy to see that 
$(\pi^*\nabla)^{(0,1)}$ defines a holomorphic structure on $\pi^*V$ provided that $(V,\nabla)$ is  an SDYM field. Then $E$ is defined to be this holomorphic vector bundle on $\mathbb{P}\setminus P$ and one can see that it has a (unique) extension to $\mathbb{P}$. \\
By the split Penrose transform in split instanton backgrounds we mean the split Penrose transform of $E$, the holomorphic vector bundle on $\mathbb{P}$ associated to a split instanton $(V,\nabla)$ on $M$.  \\
From now on suppose that  $(V,\nabla)$ is a real analytic split instanton  on $M$ and $E$ is the holomorphic vector bundle on $\mathbb{P}$ associated to it.  Since $(V,\nabla)$  is real analytic, it has an extension to a holomorphic vector bundle equipped with a holomorphic connection on an open neighborhood of $M$ in $\mathbb{M}_r$.  We denote this holomorphic vector bundle and its connection by $(V_h,\nabla_h)$. From the construction of $E$, one can see that
\begin{lemma}
The (holomorphic) vector bundles $\mu^*E$ and $\nu^*V_h$ are canonically isomorphic on some open neighborhood in $\mathbb{F}_r$ containing $G$. 
\end{lemma}
This lemma together with the first part of the split Penrose transform provides a bijection 
$$\mathcal{P}_0:H^1_c(\mathbb{P}\setminus P,E(n))\to H^1(G,\widetilde{\nu^*V_h(n)})$$
where $E(n):=E\otimes \mathcal{O}(n)$ and $\nu^*V_h(n):=\nu^*V_h\otimes\mu^*\mathcal{O}(n)$. Since $\nu^*Vh$ is trivial on the fiber of $\nu$, it is easy to compute the direct images of  $\widetilde{\nu^*V_h(n)}$. More precisely, $\nu_*^n\widetilde{\nu^*V_h(n)}=\nu_*^n\mu^*\mathcal{O}(n)\otimes V_h$. Therefore, the only remaining task to finish the split Penrose transform in the presence of split instatons is computing the differential operators. As explained before, the second and third  part of the split Penrose transform are essentially the same as ones  the complex Penrose transform. Therefore, we only need to know what the differential operators are in the complex picture. Fortunately, they have been studied  in the complex and Euclidean  picture, see \cite{BE,Ea1,Hi,We1}.  These differential operators are just the ordinary ones (i.e. when there are no split instanton backgrounds) coupled with $\nabla$. In other words, the ordinary derivatives are replaced by covariant derivatives, see \cite{Hi}.
As an example we have the following 
\begin{theorem}
Suppose that $(V,\nabla)$ is a real analytic split instanton  on $M$ and $E$ is the holomorphic vector bundle on $\mathbb{P}$ associated to it. Then there is a bijection
$$\mathcal{P}:H^1_c(\mathbb{P}\setminus P,E(-2))\to\ker\left (\square_V:\Gamma^{\omega}(
M,\widetilde{V[-1]})\to
\Gamma^{\omega}(M,\widetilde{V[-3]})\right )
$$ 
Here $V[-1]:=V\otimes \varepsilon[-1]$ and $V[-3]:=V\otimes \varepsilon[-3]$. The operator $\square_V$ is the operator obtained by coupling $\nabla$ with $\square_{2,2}$, see \cite{Hi}.

\end{theorem}
The transform in the above theorem is  the X-ray transform in the split instanton background.  
We recall how the operator $\square_V$ looks like in local coordinates. There are local coordinates $(x_1,x_2,x_3,x_4)$ on $M$ where $\square_{2,2}$ is given by
$$\square_{2,2}:=\partial_1^2+\partial_2^2-\partial_3^2-\partial_4^2$$
where $\partial_{i}:=\frac{\partial}{\partial x_{i}}$, see \cite{Ar}. If $A=\sum_{i}A_{i}dx_{i}$ is the connection form of $\nabla$ then $\square_V$, in these local coordinates, is given by 
$$\square_V=(\partial_{1}+A_{1})^2+(\partial_{2}+A_{2})^2-(\partial_{3}+A_{3})^2-(\partial_{4}+A_{4})^2$$


\end{section}

\begin{section}{Relations of the Split Penrose Transform with Representation Theory}
Finally we give some  applications of the split Penrose transform in Representation Theory and discuss the possible generalizations of it. \\
The complex Penrose transform can be used to realize some unitary representations of $SU(2,2)$ on sheaf cohomology groups on subsets of $\mathbb{CP}^3$, see \cite{BE} section 10 and references therein. In the same way, we can realize representations of $SL(4,\mathbb{R})$ via the split Penrose transform. More precisely, $SL(4,\mathbb{R})$ acts on $\mathbb{CP}^3\setminus\mathbb{RP}^3$ transitively. Hence it acts on cohomology groups $H^1_c(\mathbb{CP}^3\setminus\mathbb{RP}^3,\mathcal{O}(n))$. Via the split Penrose transform,   these cohomology groups are identified with the kernel of certain differential operators on $M$. As we saw,   $H^1_c(\mathbb{CP}^3\setminus\mathbb{RP}^3,\mathcal{O}(-2))$ is identified with the kernel of the ultra-hyperbolic operator on $M$ which is the so-called minimal representation of $SL(4,\mathbb{R})\cong SO_0(3,3)$, see \cite{KO}. Therefore, the natural action of 
$SL(4,\mathbb{R})$ on  $H^1_c(\mathbb{CP}^3\setminus\mathbb{RP}^3,\mathcal{O}(-2))$ gives a realization of the minimal representation of $SL(4,\mathbb{R})$.  Since this representation is unitarizable, it would be interesting to find an  $SL(4,\mathbb{R})$-invariant inner product on $H^1_c(\mathbb{CP}^3\setminus\mathbb{RP}^3,\mathcal{O}(-2))$. Clearly one can realize other representations of $SL(4,\mathbb{R})$ via the split Penrose transform and it would be interesting to see how these cohomological interpretations can shed light on the Representation Theory of $SL(4,\mathbb{R})$.\\
 As for the generalizations of the split Penrose transform, one can see that it easily generalizes to $SL(n,\mathbb{R})$ ($n\geq 4$). More precisely, the split Penrose transform gives a transform from compact cohomological data on $\mathbb{CP}^n\setminus\mathbb{RP}^n$ to real analytic data on $Gr(2,\mathbb{R}^n)$, the Grassmannian of 2-planes in $\mathbb{R}^n$, see \cite{BE1}.  It seems possible to derive a version of the split Penrose transform for the Funk transform (i.e. in the case of $SL(3,\mathbb{R})$) as well, see \cite{BEGM}.  \\
 More generally, one can ask for the split Penrose transform for any real semisimple (or reductive) group $G$. More precisely, the split Penrose transform must relate cohomology groups with supports on open $G$-invariant subsets of \mbox{$G_{\mathbb{C}}/P_{\mathbb{C}}\setminus G/P$} and real analytic data on $G/Q$ where $P$ and $Q$ are parabolic subgroups of $G$ and $G_{\mathbb{C}}$ and $P_{\mathbb{C}}$ are complexifications of $G$ and $P$  (we assume that $P=G\cap P_{\mathbb{C}}$), see \cite{BE}  for the generalizations of the complex Penrose transform.\\
 
 There is a transform in Representation theory due to W. Schmid which is similar to the Penrose transform, see \cite{Sc}.  
One can study the relationship of the split Penrose transform  and Schmid's transform.   The setting for  Schmid's transform is 
not exactly as the one for the split Penrose transform but nevertheless there are some similarities which indicate that there might be a close relation between them.\\

\end{section}

\end{document}